\documentclass{article}%
\usepackage{amsmath}
\usepackage{amsfonts}
\usepackage{amssymb}
\usepackage{graphicx}%
\numberwithin{equation}{section}
\newtheorem{theorem}{Theorem}[section]

\newtheorem{definition}[theorem]{Definition}

\newtheorem{lemma}[theorem]{Lemma}

\newtheorem{remark}[theorem]{Remark}

\newenvironment{proof}[1][Proof]{\noindent\textbf{#1.} }{\ \rule{0.5em}{0.5em}}
\begin{document}

\title{Estimates in Generalized Morrey Spaces for Weak Solutions to Divergence
Degenerate Parabolic Systems \thanks{This work was supported by the National
Natural Science Foundation of China (Grant Nos. 10871157 and 11001221) and
Specialized Research Found for the Doctoral Program of Higher Education (No.
200806990032).}}
\author{\ Yan DONG, Maochun ZHU, Pengcheng NIU \thanks{Corresponding author, E-mail
addresses: pengchengniu@nwpu.edu.cn}}
\maketitle

\textbf{Abstract:} Let $\mathrm{{X}}=(X_{1},\cdots,X_{q})$ be a family of real
smooth vector fields satisfying H\"{o}mander's condition. The purpose of this
paper is to establish gradient estimates in generalized Morrey spaces for
weak solutions of the divergence degenerate parabolic system related to $X$ :%
\[
u_{t}^{i}+X_{\alpha}^{\ast}(a_{ij}^{\alpha\beta}(z)X_{\beta}u^{j}%
)=g_{i}+X_{\alpha}^{\ast}f_{i}^{\alpha}(z),
\]
where $\alpha,\beta=1,2,\cdots,q,$ $i,j=1,2,\cdots,N$, $X_{\alpha}^{\ast}$ is
the transposed vector field of $X_{\alpha}$, $z=(t,x)\in{\mathbb{R}}^{n+1}$,
and coefficients $a_{ij}^{\alpha\beta}(z)$ belong to the space $VMO$ induced
by the vector fields $X_{1}, ... ,X_{q}$.


\section{Introduction}

Let us consider a family of real smooth vector fields
\[
X_{j}=\sum\limits_{k=1}^{n}{b_{jk}(x)\frac{\partial}{{\partial x_{k}}}%
},j=1,2,\cdots,q,(q\leq n),
\]
defined in a neighborhood $\tilde{\Omega}$ of some bounded domain
$\Omega\subset{\mathbb{R}}^{n}$, satisfying H\"{o}rmander's condition, namely,
the Lie algebra generated by the family $\mathrm{{X}}=(X_{1},\cdots,X_{q})$ at
any point of $\tilde{\Omega}$ spans ${{\mathbb{R}}^{n}}$, see [10].

Equations and systems involving vector fields have received much attention
during the recent years, see [1, 8, 12, 14, 15, 16, 19, 20, 21] etc.. The
Morrey estimates for elliptic systems in Carnot-Carath\'{e}odory space have
been studied by G. Di Fazio and M. Fanciullo in [6]. The aim of this paper is
to establish gradient estimates in generalized Morrey spaces for weak
solutions of the divergence degenerate parabolic system related to X. Of
course, our work is greatly influenced by those in the classic Euclidean case,
that is $X_{i}=\frac{\partial}{{\partial x_{i}}}, i=1,...,n$, where fairly
complete results have been obtained, for example, see [2, 9, 11, 18] etc.. For
parabolic system with constant coefficients, Schauder and $L^{p}$ estimates
were studied by Schlag in [17], while when coefficients are discontinuous and
belong some $VMO$ space, Mcbride in [13] derived the generalized Morrey
estimates for gradients of weak solutions. For some earlier studies, we quote
[3, 4, 18] and the references therein.

In this paper, the degenerate parabolic system we considered is of the type
\begin{equation}
u_{t}^{i}+X_{\alpha}^{\ast}(a_{ij}^{\alpha\beta}(z)X_{\beta}u^{j}%
)=g_{i}+X_{\alpha}^{\ast}f_{i}^{\alpha}(z), \label{1.1}%
\end{equation}
where $\alpha,\beta=1,2,\cdots,q;$ $i,j=1,2,\cdots,N$, $z=(x,t)\in
{{\mathbb{R}}^{n+1}}$, $X_{j}^{\ast}$ is the transposed vector field of
$X_{j}$, $X_{j}^{\ast}=-X_{j}+c_{j}$ ($c_{j}=-\sum\limits_{k=1}^{n}%
{\frac{{\partial b_{jk}}}{{\partial x_{k}}}}\in C^{\infty}(\Omega)$), $\Omega$
is a bounded domain in ${{\mathbb{R}}^{n}}$.

The main difficulty in our setting is that the presence of commutators of
vector fields which does not allow us to differentiate the equation. In order
to overcome this and apply the method in [13] to our system, we need to resort
to some conclusions proved by Xu in [20], and prove that some relative results in
the classic Euclidean case are still hold in our setting.

Our basic assumption is:

(H) Let $g_{i}$ and $f_{i}^{\alpha}$ in (\ref{1.1}) belong to the generalized
Morrey space $L_{\varphi}^{2,\lambda}(Q_{T})$, $0\leq\lambda<Q+2$ (the number
$Q$ is the homogeneous dimension relative to $\Omega$), and coefficients
$a_{ij}^{\alpha\beta}(z)$ belong to $L^{\infty}(Q_{T})\cap VMO(Q_{T})$, where
we refer to Section 2 for the precise meaning of $L_{\varphi}^{2,\lambda
}(Q_{T}),Q_{T},\varphi$ and $VMO(Q_{T})$. Also let $a_{ij}^{\alpha\beta}(x,t)$
satisfy the uniform ellipticity condition:%
\begin{equation}
\Lambda^{-1}|\xi|^{2}\leq a_{ij}^{\alpha\beta}(x,t)\xi_{\alpha}^{i}\xi_{\beta
}^{j}\leq\Lambda|\xi|^{2}, \label{1.2}%
\end{equation}
where $\Lambda>1$, $\xi\in{\mathbb{R}}^{(q+1)N}$, $(x,t)\in Q_{T}$.

We say $u\in V_{2}(Q_{T})$ (see Section 2) is a weak solution of
(\ref{1.1}), if for any vector-valued function $\psi \in
C_{0}^{\infty}(Q_{T})$,
\[
\iint\nolimits_{Q_{T}}{[u_{t}^{i}\psi+a_{ij}^{\alpha\beta}X_{\alpha}%
\psi^{i}X_{\beta}u^{j}]dz}=\iint\nolimits_{Q_{T}}{\left[  {g_{i}%
\psi+f_{i}^{\alpha}X_{\alpha}\psi}\right]  dz}.
\]

Now, we state the main result of this paper.

\begin{theorem}
\label{T1.1}Under the assumption (H), let $u\in V_{2}(Q_{T})$ be a weak
solution of (\ref{1.1}) in $Q_{T}$. Suppose that there exists $\gamma$, such
that $\lambda<\gamma<Q+2$ and the function $\frac{{r^{\gamma-\lambda}}%
}{{\varphi^{2}(r)}}$ ($r>0$) is almost increasing (see Section 2 below). Then
$Xu\in L_{\varphi}^{2,\lambda}(Q^{\prime})$ for any $Q^{\prime}\Subset
Q^{\prime\prime}\Subset Q_{T}$. Moreover, the following estimate holds
\[
\Vert{Xu}\Vert_{L_{\varphi}^{2,\lambda}(Q^{\prime})}^{2}\leq c(\Vert{Xu}%
\Vert_{L^{2}(Q^{\prime\prime})}+\Vert f\Vert_{L_{\varphi}^{2,\lambda}(Q_{T}%
)}+\Vert g\Vert_{L_{\varphi}^{2,\lambda}(Q_{T})}).
\]

\end{theorem}

The plan of the paper is organized as follows. In Section 2, we introduce some
function spaces such as generalized Morrey spaces, generalized Sobolev spaces,
and give some known results which will be used. Section 3 is devoted to
deducing a Caccioppoli inequality (Lemma 3.1) and $L^{2}$ estimates for
derivatives (with respect to vector fields and $t$) of weak solutions of (3.1)
(Lemma 3.4). Using the reverse H\"{o}lder inequality on the homogeneous space,
we prove a higher integrability to (3.1)(see Theorem 3.9). With the
help of the results in Sections 2 and 3, we complete the proof of Theorem 1.1 in
Section 4.

\section{Preliminaries}

In this section we introduce some preparatory material related to
H\"{o}rmander's vector fields and state some function spaces. Several known
results which will be used later are collected.

For every multi-index $I=\left(  {i_{1},i_{2},\ldots,i_{k}}\right) $, we
denote the length of $I$ by $\left\vert I\right\vert =k$, and set
$\ \ \ \ \ \ \ \ \ \ \ \ $%
\[
X_{I}=X_{i_{i}}X_{i_{2}}\ldots X_{i_{k}},X_{\beta}=[X_{\beta_{d}},
[X_{\beta_{d-1}},\cdots\lbrack X_{\beta_{2}},X_{\beta_{1}}]\cdots]].
\]
The length of commutator $X_{\beta}$ is denoted by $\left\vert \beta
\right\vert =d$.

\begin{definition}
\label{D2.2}(Carnot-Carath\'{e}odory distance). An absolutely continuous curve
$\gamma:[0,T]\rightarrow\Omega$ is called a sub-unit curve with respect to the
system $\mathrm{{X}}$, if $\gamma^{\prime}(t)$ exists and
satisfies that for any $\xi\in{\mathbb{R}}^{n},$
\[
<\gamma^{\prime}(t),\xi>^{2}\leq\sum\limits_{j=1}^{q}{<X_{j}(\gamma
(t)),\xi>^{2}},\text{{a.e. }}t\in\lbrack0,T].
\]
The length of $\gamma$ is denoted by $l_{S}\left(  \gamma\right)  =T$. Given
any $x,y\in\Omega$, we stand for the collection of all sub-unit curves
connecting $x$ and $y$ by $\Phi(x,y)$ and define
\[
d_{\mathrm{{X}}}(x,y)=\inf\{l_{S}(\gamma):\gamma\in\Phi(x,y)\}.
\]
Note that the function $d_{\mathrm{{X}}}(x,y)$ is finite for any $x,y\in
\Omega$, and $d_{\mathrm{{X}}}$ is really a distance in $\Omega$. One calls
that $d_{\mathrm{{X}}}$ is a Carnot-Carath\'{e}odory distance.
\end{definition}

A metric ball of center $x$ and radius $R$ is denoted by
\[
B_{R}{(x)}=B(x,R)=\{y\in\Omega:d_{X}(x,y)<R\}.
\]
When we do not consider the center of a ball, we will simply write $B_{R}$
instead of $B(x,R)$.

Due to \cite{b14}, for $\Omega\subset{\mathbb{R}}^{n}$, there exist constants
$C_{D},R_{D}>0$ such that for every $x_{0}\in\Omega$ and $0<R<R_{D}$, one has%
\[
|{B(x_{0},2R)}|\leq C_{D}|{B(x_{0},R)}|.
\]
Moreover, for every $R\leq R_{D}$ and $\tau\in(0,1)$, we have
\begin{equation}
\left\vert {B_{\tau R}}\right\vert \geq C_{_{D}}^{-1}\tau^{Q}\left\vert
{B_{R}}\right\vert . \label{2.1}%
\end{equation}
Through out of this paper, we denote $Q_{T}=\Omega\times(0,T]$ and
$z_{0}=(x_{0},t_{0})\in Q_{T}\subset{\mathbb{R}}^{n+1}$ . A parabolic cylinder
with vertex at $z_{0}$ is denoted by
\[
Q_{R}(z_{0})=B_{R}(x_{0})\times(t_{0}-R^{2},t_{0}].
\]

In the sequel, let us denote $I_{R}\left(  {t_{0}}\right)  =(t_{0}-R^{2}%
,t_{0}]$ and the parabolic boundary of $Q_{R}$\ by $\partial_{p}Q_{R}$.
Denote the Lebesgue measure of $B(x,R)$ in the $n$-dimensional space by
$\left\vert {B(x,R)}\right\vert $, and the Lebesgue measure of $Q_{R}(z_{0})$
in the $n+1$-dimensional space by $\left\vert {Q_{R}(z_{0})}\right\vert $.

\begin{definition}
\label{D2.3}(Almost increasing function, see \cite{b11}). A function
$h:[0,d_{0}]\rightarrow\lbrack0,\infty)$ is said almost increasing, where
$d_{0}>0$, if there exists $K_{h}\geq1$, such that for any $0\leq s\leq t\leq
d_{0}$, the following holds
\[
h(s)\leq K_{h}h(t).
\]

\end{definition}

\begin{definition}
\label{D2.4} Let $1\leq p<+\infty, 0\leq\lambda<Q+2$ and $\varphi$ be a
continuous function on $[0,d]$ such that $\varphi>0$ on $(0,d]$, where $d$ is
the diameter of $Q_{T}$. We say that $f\in L^{p}(Q_{T})$ belongs to a
generalized Morrey space $L_{\varphi}^{p,\lambda}(Q_{T})$, if
\[
\left\Vert f\right\Vert _{L_{\varphi}^{p,\lambda}}=\mathop {\sup
}\limits_{z_{0}\in Q_{T},0\leq\rho\leq d}\frac{1}{{\varphi(\rho)}}%
(\rho^{-\lambda}\iint\nolimits_{Q_{T}\cap Q_{\rho}(z_{0})}{\left\vert
f\right\vert ^{p}dz})^{\frac{1}{p}}<\infty.
\]

\end{definition}

It is easy to prove that the space $L_{\varphi}^{p,\lambda}(Q_{T})$ is a Banach space
as in \cite{b13}.

\begin{definition}
\label{D2.5}(BMO and VMO space). For any $f\in L^{1}(Q_{T})$, we set
\[
\eta\left(  r\right)  =\mathop {\sup }\limits_{z_{0}\in Q,0\leq\rho\leq
r}(\frac{1}{{\left\vert {Q_{T}\cap Q_{\rho}(z_{0})}\right\vert }}%
\iint\nolimits_{Q_{T}\cap Q_{\rho}(z_{0})}{\left\vert {f(z)-f_{Q_{T}\cap
Q_{\rho}(z_{0})}(z)}\right\vert dz}),
\]
where $f_{Q_{T}\cap Q_{\rho}(z_{0})}=\frac{1}{{\left\vert {Q_{T}\cap
Q_{\rho }(z_{0})}\right\vert }}\iint\nolimits_{Q_{T}\cap
Q_{\rho}(z_{0})}{f(z)dz}$. If $\mathop {\sup
}\limits_{r>0}\eta\left(  r\right)  <\infty$, then it says $f\in BMO(Q_{T}%
)$(Bounded Mean Oscillation). Moreover, if $\eta\left(  r\right)
\rightarrow0$ as $r\rightarrow0$, then we call $f\in VMO(Q_{T})$(Vanishing
Mean Oscillation).
\end{definition}

\begin{definition}
\label{D2.6}(Generalized Sobolev space). The space
\[
V_{2}(Q_{T})=\{u:u\in L^{\infty}(0,T;L^{2}(Q_{T})),Xu\in L^{2}(Q_{T})\}
\]
is called a generalized Sobolev space, where $\left\vert
Xu\right\vert =(\sum\limits_{i=1}^{q}{\left\vert {X_{i}u}\right\vert
^{2}})^{\frac{1}{2}}$.
\end{definition}

\begin{lemma}
\label{l2.1}(see \cite{b11}). Let $H$ be a non-negative almost increasing
function in $[0,R_{0}]$ and $F$ a positive function on $(0,R_{0}]$. Suppose
that $H$ and $F$ satisfy

(1) There exist positive constants $A,B,\varepsilon$ and $\beta$ such that for any
$0\leq\rho\leq R\leq R_{0}$,
\begin{equation}
H(\rho)\leq(A(\frac{\rho}{R})^{\beta}+\varepsilon)H(R)+BF(R); \label{2.2}%
\end{equation}
\

(2) There exists $\gamma\in(0,\beta)$ such that $\frac{{\rho^{\gamma}}%
}{{F(\rho)}}$ is almost increasing in $(0,R_{0}]$.

Then there exist $\varepsilon_{0}=\varepsilon_{0}(A,\beta,\gamma)$ and
$C=C(A,\beta,\gamma,K_{H})$ such that if $\varepsilon<\varepsilon_{0}$, one
has
\begin{equation}
H(\rho)\leq C\frac{{F(\rho)}}{{F(R)}}H(R)+CBF(\rho). \label{2.3}%
\end{equation}

\end{lemma}

The following technical lemma is from \cite{b5}.

\begin{lemma}
\label{l2.2} Let $f(t)$ be a bounded nonnegative function on $[T_{0},T_{1}]$, $T_{1}>T_{0}\geq0$. Suppose that for any $s$ and $t$, $T_{0}\leq t<s\leq
T_{1}$, $f$ satisfies
\[
f(t)\leq\theta f(s)+\frac{A}{{(s-t)^{\alpha}}}+B,
\]
where $\theta,A,B,\alpha$ are nonnegative constants and $\theta<1$. Then for
any $T_{0}\leq\rho<R\leq T_{1}$, one has
\[
f(\rho)\leq C[\frac{A}{{(R-\rho)^{\alpha}}}+B],
\]
where $C$ depends only on $\alpha$.
\end{lemma}

\section{Homogeneous parabolic system with constant coefficients}

Let us consider the homogeneous degenerate parabolic system
\begin{equation}
u_{t}^{i}+X_{\alpha}^{\ast}(a_{ij}^{\alpha\beta}X_{\beta}u^{j})=0, \label{3.1}%
\end{equation}
where coefficients $a_{ij}^{\alpha\beta}$ are constants and satisfy
(\ref{1.2}). We will establish a Caccioppoli inequality and $L^{2}$ estimates for
derivatives (with respect to vector fields $X_{1},...,X_{q}$ and the variable
$t$) of weak solutions of (\ref{3.1}) by extending results in [20]. Using the
reverse H\"{o}der inequality on the homogeneous space, a higher integrability
to (\ref{3.1}) is proved. To simplify the notations, in the sequel,
$Q_{R}\left(  {z_{0}}\right)  $, $B_{R}\left(  {x_{0}}\right)  $,
$I_{R}\left(  {t_{0}}\right) $ and $dxdt$ are written as $Q_{R}$, $B_{R}$,
$I_{R}$ and $dz$, respectively.

\begin{lemma}
\label{l3.1}(Caccioppoli inequality). Let $u\in V_{2}(Q_{T})$ be a weak
solution of (\ref{3.1}). Then for any $Q_{R}\subset Q_{T}$ and $\rho<R$,
\begin{equation}
\mathop {\sup }\limits_{I_{\rho}}\int_{B_{\rho}}{\left\vert u\right\vert ^{2}%
}dx+\iint\nolimits_{Q_{\rho}}{\left\vert {Xu}\right\vert ^{2}dxdt}\leq\frac
{c}{{\left(  {R-\rho}\right)  ^{2}}}\iint\nolimits_{Q_{R}}{\left\vert
u\right\vert ^{2}dxdt}. \label{3.2}%
\end{equation}
Furthermore, for any $b\in{\mathbb{R}}$, it follows
\begin{equation}
\mathop {\sup }\limits_{I_{\rho}}\int_{B_{\rho}}{\left\vert {u-b}\right\vert
^{2}}dx+\iint\nolimits_{Q_{\rho}}{\left\vert {Xu}\right\vert ^{2}dxdt}%
\leq\frac{c}{{\left(  {R-\rho}\right)  ^{2}}}\iint\nolimits_{Q_{R}}{\left\vert
{u-b}\right\vert ^{2}dxdt} \label{3.3}%
\end{equation}
and
\begin{equation}
\iint\nolimits_{Q_{\rho}}{|{u_{t}}|^{2}dxdt}\leq\frac{c}{{({R-\rho})^{4}}%
}\iint\nolimits_{Q_{R}}{|u|^{2}dxdt}. \label{3.4}%
\end{equation}

\end{lemma}

\begin{proof}
Given $B_{\rho}\subset B_{R}\subset\Omega$, choose a test function ${f_{i}%
}(x)=u^{i}\xi^{2}(x)\eta(t)$ with
\[
\xi(x)\in C_{0}^{\infty}(B_{R}), 0\leq\xi\leq1,\left\vert {X\xi}\right\vert
\leq\frac{C}{{R-\rho}}, \xi=1(\text{in }B_{\rho})
\]
and
\[
\eta(t)=\left\{  {\begin{array}{*{20}c}
{\frac{{t - \left( {t_0  - R^2 } \right)}}{{R^2  - \rho ^2 }},t \in (t_0  - R^2 ,t_0  - \rho ^2 ),}  \\
{1,{\kern 1pt} {\kern 1pt} {\kern 1pt} {\kern 1pt} {\kern 1pt} {\kern 1pt} {\kern 1pt} {\kern 1pt} {\kern 1pt} {\kern 1pt} {\kern 1pt} {\kern 1pt} {\kern 1pt} {\kern 1pt} {\kern 1pt} {\kern 1pt} {\kern 1pt} {\kern 1pt} {\kern 1pt} {\kern 1pt} {\kern 1pt} {\kern 1pt} {\kern 1pt} {\kern 1pt} {\kern 1pt} {\kern 1pt} {\kern 1pt} {\kern 1pt} {\kern 1pt} {\kern 1pt} {\kern 1pt} t \in [t_0  - \rho ^2 ,t_0 ).}  \\
\end{array}}\right.
\]
Multiplying both sides of (\ref{3.1}) by ${f_{i}}(x)$ and integrating on
$Q_{R}^{\prime}=B_{R}(x_{0})\times(t_{0}-R^{2},t]$, we get%
\begin{align*}
0  &  =\iint\nolimits_{Q_{R}^{\prime}}{[u_{t}^{i}+X_{\alpha}^{\ast}%
(a_{ij}^{\alpha\beta}X_{\beta}u^{j})]u^{i}\xi^{2}\eta dz}\\
&  =\iint\nolimits_{Q_{R}^{\prime}}{\left[  {(\frac{1}{2}\left\vert
u\right\vert ^{2}\eta)_{t}\xi^{2}-\frac{1}{2}\left\vert u\right\vert ^{2}%
\xi^{2}\eta_{t}+a_{ij}^{\alpha\beta}\xi^{2}\eta X_{\alpha}u^{i}X_{\beta}u^{j}%
}\right.  }\\
&  \left.  {+2a_{ij}^{\alpha\beta}u^{i}\xi\eta X_{\alpha}\xi X_{\beta}u^{j}%
}\right]  dz,
\end{align*}
and then by moving terms,
\begin{align}
&  \iint\nolimits_{Q_{R}^{\prime}}{[(\frac{1}{2}\left\vert u\right\vert
^{2}\eta)_{t}\xi^{2}+a_{ij}^{\alpha\beta}\xi^{2}\eta X_{\alpha}u^{i}X_{\beta
}u^{j}]dz}\nonumber\\
&  =\iint\nolimits_{Q_{R}^{\prime}}{[\frac{1}{2}\left\vert u\right\vert
^{2}\xi^{2}\eta_{t}-2a_{ij}^{\alpha\beta}u^{i}\xi\eta X_{\alpha}\xi X_{\beta
}u^{j}]dz}\nonumber\\
&  \leq\iint\nolimits_{Q_{R}^{\prime}}{\frac{1}{2}\left\vert u\right\vert
^{2}\xi^{2}\eta_{t}dz}+\varepsilon\iint\nolimits_{Q_{R}^{\prime}}{\xi^{2}%
\eta\left\vert {Xu}\right\vert ^{2}dz}+C_{\varepsilon}\iint\nolimits_{Q_{R}%
^{\prime}}{\eta\left\vert u\right\vert ^{2}\left\vert {X\xi}\right\vert
^{2}dz}. \label{3.5}%
\end{align}
Using
\begin{align*}
\iint\nolimits_{Q_{R}^{\prime}}{(\frac{1}{2}\left\vert u\right\vert ^{2}%
\eta)_{t}\xi^{2}dz}  &  =\int_{B_{R}}{\int_{(t_{0}-R^{2},t]}{(\frac{1}%
{2}\left\vert u\right\vert ^{2}\eta)_{t}\xi^{2}dt}dx}\\
&  =\eta\left(  t\right)  \int_{B_{R}}{\frac{1}{2}\left\vert u\right\vert
^{2}\xi^{2}dx},
\end{align*}
and (\ref{1.2}), we have from (\ref{3.5}) that%
\begin{align*}
&  \eta\left(  t\right)  \int_{B_{R}}{\frac{1}{2}\left\vert u\right\vert
^{2}\xi^{2}dx}+C\iint\nolimits_{Q_{R}^{\prime}}{\xi^{2}\eta\left\vert
{Xu}\right\vert ^{2}dz}\\
&  \leq\iint\nolimits_{Q_{R}^{\prime}}{\frac{1}{2}\left\vert u\right\vert
^{2}\xi^{2}\eta_{t}dz}+\varepsilon\iint\nolimits_{Q_{R}^{\prime}}{\xi^{2}%
\eta\left\vert {Xu}\right\vert ^{2}dz}+C_{\varepsilon}\iint\nolimits_{Q_{R}%
^{\prime}}{\eta\left\vert u\right\vert ^{2}\left\vert {X\xi}\right\vert
^{2}dz}.
\end{align*}
In the light of properties of $\xi\left(  x\right) $ and $\eta\left(
t\right)  $, it implies%
\begin{align*}
&  \eta\left(  t\right)  \int_{B_{R}}{\left\vert u\right\vert ^{2}\xi^{2}%
dx}+\iint\nolimits_{Q_{R}^{\prime}}{\xi^{2}\eta\left\vert {Xu}\right\vert
^{2}dz}\\
&  \leq C_{\varepsilon}\iint\nolimits_{Q_{R}^{\prime}}{\left\vert u\right\vert
^{2}\xi^{2}\eta_{t}dz}+C_{\varepsilon}\iint\nolimits_{Q_{R}^{\prime}}%
{\eta\left\vert u\right\vert ^{2}\left\vert {X\xi}\right\vert ^{2}dz}\\
&  \leq C_{\varepsilon}\iint\nolimits_{Q_{R}^{\prime}}{\left\vert u\right\vert
^{2}(\frac{1}{{R^{2}-\rho^{2}}}+\frac{C}{{(R-\rho)^{2}}})dz},
\end{align*}
thus
\[
\mathop {\sup }\limits_{I_{\rho}}\int_{B_{\rho}}{\left\vert u\right\vert ^{2}%
}dx+\iint\nolimits_{Q_{\rho}}{\left\vert {Xu}\right\vert ^{2}dxdt}\leq\frac
{c}{{\left(  {R-\rho}\right)  ^{2}}}\iint\nolimits_{Q_{R}}{\left\vert
u\right\vert ^{2}dxdt},
\]
namely, (\ref{3.2}) is proved.

The proof of (\ref{3.3}) is similar to that of (\ref{3.2}), just taking the
test function ${f_{i}}(x)=(u^{i}-b)\xi^{2}(x)\eta(t)$ instead. We omit the details.

Now we come to prove (\ref{3.4}). Let $\rho\leq s<l\leq R$ with $l-s=s-\rho$
and ${f_{i}}(x)=u_{t}^{i}\xi_{1}^{2}(x)\eta_{1}(t)$ be a test function with
\[
\xi_{1}(x)\in C_{0}^{\infty}(B_{s}), 0\leq\xi_{1}\leq1,\left\vert {X\xi_{1}%
}\right\vert \leq\frac{C}{{s-\rho}}, \xi_{1}=1 (\text{in }B_{\rho})
\]
and
\[
\eta_{1}(t)=\left\{  {\begin{array}{*{20}c}
{\frac{{t - \left( {t_0  - s^2 } \right)}}{{s^2  - \rho ^2 }},t \in (t_0  - s^2 ,t_0  - \rho ^2 )},  \\
{1,{\kern 1pt} {\kern 1pt} {\kern 1pt} {\kern 1pt} {\kern 1pt} {\kern 1pt} {\kern 1pt} {\kern 1pt}
{\kern 1pt} {\kern 1pt} {\kern 1pt} {\kern 1pt} {\kern 1pt} {\kern 1pt} {\kern 1pt} {\kern 1pt}
{\kern 1pt} {\kern 1pt} {\kern 1pt} {\kern 1pt} {\kern 1pt} {\kern 1pt} {\kern 1pt} {\kern 1pt}
{\kern 1pt} {\kern 1pt} {\kern 1pt} {\kern 1pt} {\kern 1pt} {\kern 1pt} {\kern 1pt} t \in [t_0  - \rho ^2 ,t_0 )}.  \\
\end{array}}\right.
\]
Multiplying both sides in (\ref{3.1}) by ${f_{i}}(x)$ and integrating on $Q_{s}$,
one gets%
\begin{align*}
0  &  =\iint\nolimits_{Q_{s}}{[(u_{t}^{i})^{2}\xi_{1}^{2}\eta_{1}%
+a_{ij}^{\alpha\beta}\xi_{1}^{2}\eta_{1}X_{\alpha}u_{t}^{i}X_{\beta}%
u^{j}+2a_{ij}^{\alpha\beta}u_{t}^{i}\xi_{1}\eta_{1}X_{\alpha}\xi_{1}X_{\beta
}u^{j}]dz},
\end{align*}
then%
\begin{align*}
\iint\nolimits_{Q_{s}}{\xi_{1}^{2}\eta_{1}\left\vert {u_{t}}\right\vert
^{2}dz}  &  \leq C\iint\nolimits_{Q_{s}}{\left[  {\xi_{1}^{2}\eta
_{1}\left\vert {Xu_{t}}\right\vert \left\vert {Xu}\right\vert +\eta
_{1}\left\vert {\xi_{1}u_{t}}\right\vert \left\vert {Xu}\right\vert \left\vert
{X\xi_{1}}\right\vert }\right]  dz}\\
&  \leq{\varepsilon} \iint\nolimits_{Q_{s}}{\xi_{1}^{2}\eta_{1}\left\vert
{Xu_{t}}\right\vert ^{2}dz}+\frac{C}{{\varepsilon}}\iint\nolimits_{Q_{s}}%
{\xi_{1}^{2}\eta_{1}\left\vert {Xu}\right\vert ^{2}dz}\\
&  + \frac{1}{2}\iint\nolimits_{Q_{s}}{\eta_{1}\left\vert {\xi_{1}u_{t}%
}\right\vert ^{2}dz} + C \iint\nolimits_{Q_{s}}{\eta_{1}\left\vert
{Xu}\right\vert ^{2}\left\vert {X\xi_{1}}\right\vert ^{2}dz}.
\end{align*}
Noting properties of $\xi_{1}\left(  x\right)  $ and $\eta_{1}\left(
t\right)  $, it yields%
\begin{align}
\iint\nolimits_{Q_{\rho}}{\left\vert {u_{t}}\right\vert ^{2}dz}  &
\leq2{\varepsilon} \iint\nolimits_{Q_{s}}{\left\vert {Xu_{t}}\right\vert
^{2}dz}\nonumber\\
&  +\frac{C}{{\varepsilon}}\iint\nolimits_{Q_{s}}{\left\vert {Xu}\right\vert
^{2}dz}+\frac{C}{{(s-\rho)^{2}}}\iint\nolimits_{Q_{s}}{\left\vert
{Xu}\right\vert ^{2}dz}. \label{3.6}%
\end{align}
Since $u_{t}$ is still a weak solution of (\ref{3.1}), we apply (\ref{3.2}) to
$u_{t}$ and have
\[
\iint\nolimits_{Q_{s}}{\left\vert {Xu_{t}}\right\vert ^{2}dz}\leq\frac
{C}{{(l-s)^{2}}}\iint\nolimits_{Q_{l}}{\left\vert {u_{t}}\right\vert ^{2}dz}%
\]
and
\[
\iint\nolimits_{Q_{s}}{\left\vert {Xu}\right\vert ^{2}dz}\leq\frac
{C}{{(l-s)^{2}}}\iint\nolimits_{Q_{l}}{\left\vert u\right\vert ^{2}dz}.
\]
Inserting the above two inequalities into (\ref{3.6}) and using $l-s=s-\rho$,
it obtains%
\begin{align*}
&  \iint\nolimits_{Q_{\rho}}{\left\vert {u_{t}}\right\vert ^{2}dz}\\
&  \leq\frac{{2{\varepsilon} C}}{{(l-s)^{2}}}\iint\nolimits_{Q_{l}}{\left\vert
{u_{t}}\right\vert ^{2}dz}+\frac{C}{{{\varepsilon} (l-s)^{2}}}\iint
\nolimits_{Q_{l}}{\left\vert u\right\vert ^{2}dz}+\frac{C}{{(s-\rho
)^{2}(l-s)^{2}}}\iint\nolimits_{Q_{l}}{\left\vert u\right\vert ^{2}dz}\\
& \leq\frac{{{\varepsilon} C}}{{(l-s)^{2}}}\iint\nolimits_{Q_{l}}{\left\vert
{u_{t}}\right\vert ^{2}dz}+\frac{C}{{{\varepsilon}(l-s)^{2}}}\iint
\nolimits_{Q_{l}}{\left\vert u\right\vert ^{2}dz}+\frac{C}{{(l-s)^{4}}}%
\iint\nolimits_{Q_{l}}{\left\vert u\right\vert ^{2}dz}.
\end{align*}
Taking $\varepsilon=\frac{{(l-s)^{2}}}{{4C}}$, it follows
\[
\iint\nolimits_{Q_{\rho}}{\left\vert {u_{t}}\right\vert ^{2}dz}\leq\frac{1}%
{4}\iint\nolimits_{Q_{l}}{\left\vert {u_{t}}\right\vert ^{2}dz}+\frac
{C}{{(l-\rho)^{4}}}\iint\nolimits_{Q_{l}}{\left\vert u\right\vert ^{2}dz}%
\]
and then (3.4) from Lemma \ref{l2.2}.
\end{proof}

\begin{remark}
\label{R3.1}Checking carefully the proof of Lemma \ref{l3.1}, one find that
conclusions in Lemma \ref{l3.1} are still hold for the homogeneous parabolic
system with variable coefficients, provided coefficients are bounded and
satisfy (\ref{1.2}). It will be used in Section 4.
\end{remark}

\begin{lemma}
\label{l3.2} Let $u\in C^{\infty}(Q_{T})$, $B_{R}\subset\Omega$ and
$I_{R}\subset\left(  {0,T}\right) $. Then
\end{lemma}

(i) when $k>\frac{Q}{2}$, there exist positive constants $R_{0}$ and $c$ such
that for any $R\leq R_{0}$,
\begin{equation}
\mathop {\sup }\limits_{x\in B_{R/4}}\left\vert {u\left(  {x,t}\right)
}\right\vert \leq c\left\vert {B_{R}}\right\vert ^{-\frac{1}{2}}%
\sum\limits_{\left\vert I\right\vert \leq k}{R^{\left\vert I\right\vert
}\left\Vert {X_{I}u\left(  {x,t}\right)  }\right\Vert _{L^{2}(B_{R})}}.
\label{3.7}%
\end{equation}

(ii) when $k>1$, there exist positive constants $R_{0}$ and $c$ such that for
any $R\leq R_{0}$,
\begin{equation}
\mathop {\sup }\limits_{t\in I_{R/4}}\left\vert {u\left(  {x,t}\right)
}\right\vert \leq c\sum\limits_{2m\leq k}{R^{2m-1}\left\Vert {\partial_{_{t}%
}^{m}u\left(  {x,t}\right)  }\right\Vert _{L^{2}(I_{R})}}. \label{3.8}%
\end{equation}

The first statement is from Proposition 2.4 in [20]. The second is easily
proved by the same way in [20]. We omit it here.

\begin{lemma}
\label{l3.3}Let $u\in V_{2}(Q_{T})$ be a weak solution of (\ref{3.1}). Then
$u\in C^{\infty}(Q_{T})$ and for any positive integer $k$, it follows
\begin{equation}
\sum\limits_{\left\vert I\right\vert \leq k}\iint\nolimits_{Q_{R/2^{k}}%
}{{\left\vert {X_{I}u}\right\vert ^{2}dz}}\leq\frac{c}{{R^{2k}}}%
\iint\nolimits_{Q_{R}}{\left\vert u\right\vert ^{2}dz} \label{3.9}%
\end{equation}
and
\begin{equation}
\sum\limits_{\left\vert I\right\vert +2m\leq k}\iint\nolimits_{Q_{R/2^{k}}%
}{{\left\vert {X_{I}\partial_{t}^{m}u}\right\vert ^{2}dz}}\leq\frac{c}%
{{R^{2k}}}\iint\nolimits_{Q_{R}}{\left\vert u\right\vert ^{2}dz}. \label{3.10}%
\end{equation}

\end{lemma}

\begin{proof}
Denote $M^{k}(\Omega)=\{u\in L^{2}(\Omega),X_{I}u\in L^{2}(\Omega),\left\vert
I\right\vert \leq k\}$ and $Lu=u_{t}^{i}+X_{\alpha}^{\ast}(a_{ij}^{\alpha
\beta}X_{\beta}u^{j})$. Since $u$ is a weak solution of (\ref{3.1}) and $L$ is
hypoelliptic, we deduce that $u$ belongs to $C^{\infty}(Q_{T})$ from $Lu=0$.

Let us test (\ref{3.9}) by the induction on $k$. When $k=1$, setting
$\rho=\frac{R}{2}$ in (\ref{3.2}) leads to
\[
\iint\nolimits_{Q_{R/2}}{\left\vert {X_{I}u}\right\vert ^{2}dz}\leq\frac
{c}{{R^{2}}}\iint\nolimits_{Q_{R}}{\left\vert u\right\vert ^{2}dz}.
\]
Assuming that (\ref{3.9}) is true if $\left\vert I\right\vert \leq
k-1$($k\geq2$ ), we show that (\ref{3.9}) is still true when $\left\vert
I\right\vert =k$.

Let $\xi(x)\eta\left(  t\right)  $ be a cutoff function with
\[
\xi(x)\in C_{0}^{\infty}(B_{R/2^{k-1}}),0\leq\xi\leq1,\left\vert {X_{I}\xi
}\right\vert \leq\frac{C}{{R^{\left\vert I\right\vert }}},\xi=1(\text{in
}B_{R/2^{k}})
\]
and
\[
\eta(t)=\left\{  {\begin{array}{*{20}c}
{\frac{{t - \left( {t_0  - (R/2^{k - 1} )^2 } \right)}}{{(R/2^{k - 1} )^2  - (R/2^k )^2 }},t \in (t_0  - (R/2^{k - 1} )^2 ,t_0  - (R/2^k )^2 ),}  \\
{1,{\kern 1pt} {\kern 1pt} {\kern 1pt} {\kern 1pt} {\kern 1pt} {\kern 1pt} {\kern 1pt} {\kern 1pt} {\kern 1pt} {\kern 1pt} {\kern 1pt} {\kern 1pt} {\kern 1pt} {\kern 1pt} {\kern 1pt} {\kern 1pt} {\kern 1pt} {\kern 1pt} {\kern 1pt} {\kern 1pt} {\kern 1pt} {\kern 1pt} {\kern 1pt} {\kern 1pt} {\kern 1pt} {\kern 1pt} {\kern 1pt} {\kern 1pt} {\kern 1pt} {\kern 1pt} {\kern 1pt} t \in [t_0  - (R/2^k )^2 ,t_0 ).}  \\
\end{array}}\right.
\]
Denote $\tilde{L}u=a_{ij}^{\alpha\beta}X_{\alpha}X_{\beta}u^{j}$. Recalling
$Lu=0$ and $X_{\alpha}^{\ast}=-X_{\alpha}+c_{\alpha}$, one sees
\[
\tilde{L}u=a_{ij}^{\alpha\beta}c_{\alpha}X_{\beta}u^{j}+u_{t}^{i}.
\]
Due to regularity result by Rothschild and Stein (\cite{b16}) to the operator
$\tilde{L}$, we have%
\begin{align}
&  \sum\limits_{\left\vert I\right\vert =k}\iint\nolimits_{Q_{R/2^{k}}%
}{{\left\vert {X_{I}u}\right\vert ^{2}}dz}\nonumber\\
&  \leq\sum\limits_{\left\vert I\right\vert \leq k}{\int_{I_{R/2^{k}}}%
}\left\Vert {{{X_{I}u}}}\right\Vert {{{{}}_{L^{2}(B_{R/2^{k}})}^{2}}dt}%
\leq\int_{I_{R/2^{k-1}}}{\left\Vert {\xi\eta u}\right\Vert _{M^{k}%
(B_{R/2^{k-1}})}^{2}}dt\nonumber\\
&  \leq c\int_{I_{R/2^{k-1}}}{[\left\Vert {\tilde{L}(\xi\eta u)}\right\Vert
_{M^{k-2}(B_{R/2^{k-1}})}^{2}+\left\Vert {\xi\eta u}\right\Vert _{L^{2}%
(B_{R/2^{k-1}})}^{2}]}dt\nonumber\\
&  \leq c\int_{I_{R/2^{k-1}}}{[\left\Vert {a_{ij}^{\alpha\beta}c_{\alpha
}X_{\beta}(\xi\eta u)^{j}+(\xi\eta u)_{t}^{i}}\right\Vert _{M^{k-2}%
(B_{R/2^{k-1}})}^{2}+\left\Vert {\xi\eta u}\right\Vert _{L^{2}(B_{R/2^{k-1}}%
)}^{2}]}dt\nonumber\\
&  \leq c\int_{I_{R/2^{k-1}}}{[\left\Vert {a_{ij}^{\alpha\beta}c_{\alpha
}X_{\beta}(\xi\eta u)^{j}}\right\Vert _{M^{k-2}(B_{R/2^{k-1}})}^{2}%
+}\nonumber\\
&  +{\left\Vert {(\xi\eta u)_{t}^{i}}\right\Vert _{M^{k-2}(B_{R/2^{k-1}})}%
^{2}+\left\Vert {\xi\eta u}\right\Vert _{L^{2}(B_{R/2^{k-1}})}^{2}]}dt.
\label{3.11}%
\end{align}
Let us denote
\[
I\equiv c\int_{I_{R/2^{k-1}}}{\left\Vert {a_{ij}^{\alpha\beta}c_{\alpha
}X_{\beta}(\xi\eta u)^{j}}\right\Vert _{M^{k-2}(B_{R/2^{k-1}})}^{2}dt,}%
\]%
\[
II\equiv c\int_{I_{R/2^{k-1}}}{\left\Vert {(\xi\eta u)_{t}^{i}}\right\Vert
_{M^{k-2}(B_{R/2^{k-1}})}^{2}dt.}%
\]
We first estimate $II$. By properties of $\xi\left(  x\right)  $ and
$\eta\left(  t\right)  $,
\begin{align}
II  &  =c\int_{I_{R/2^{k-1}}}{\left\Vert {(\xi\eta u)_{t}}\right\Vert
_{M^{k-2}(B_{R/2^{k-1}})}^{2}}dt\nonumber\\
&  =c\sum\limits_{\left\vert I\right\vert \leq k-2}\iint
\nolimits_{Q_{R/2^{k-1}}}{{\left\vert {X_{I}(\xi\eta_{t}u+\xi\eta u_{t}%
)}\right\vert ^{2}}dz}\nonumber\\
&  \leq c\sum\limits_{\left\vert I\right\vert \leq k-2}\iint
\nolimits_{Q_{R/2^{k-1}}}{{[\left\vert {\eta_{t}uX_{I}\xi}\right\vert
^{2}+\left\vert {\xi\eta_{t}X_{I}u}\right\vert ^{2}+\left\vert {\eta
u_{t}X_{I}\xi}\right\vert ^{2}+\left\vert {\xi\eta X_{I}u_{t}}\right\vert
^{2}]}dz},\nonumber\\
&  \leq c\sum\limits_{\left\vert I\right\vert \leq k-2}{\frac{1}%
{{R^{2(\left\vert I\right\vert +2)}}}\iint\nolimits_{Q_{R/2^{k-1}}}{\left\vert
u\right\vert ^{2}dz}}+c\sum\limits_{\left\vert I\right\vert \leq k-2}{\frac
{1}{{R^{4}}}\iint\nolimits_{Q_{R/2^{k-1}}}{\left\vert {X_{I}u}\right\vert
^{2}dz}}\nonumber\\
&  {+c\sum\limits_{\left\vert I\right\vert \leq k-2}{\frac{1}{{R^{2\left\vert
I\right\vert }}}\iint\nolimits_{Q_{R/2^{k-1}}}{\left\vert {u_{t}}\right\vert
^{2}dz}}+c\sum\limits_{\left\vert I\right\vert \leq k-2}\iint
\nolimits_{Q_{R/2^{k-1}}}{{\left\vert {X_{I}u_{t}}\right\vert ^{2}dz}}.}
\label{3.12}%
\end{align}
From the assertion for $\left\vert I\right\vert \leq k-1$, it follows
\begin{equation}
\sum\limits_{\left\vert I\right\vert \leq k-2}\iint\nolimits_{Q_{R/2^{k-1}}%
}{{\left\vert {X_{I}u}\right\vert ^{2}dz}}\leq\frac{c}{{R^{2(k-2)}}}%
\iint\nolimits_{Q_{R}}{\left\vert u\right\vert ^{2}dz}. \label{3.13}%
\end{equation}
On the other hand, since $u_{t}$ is still a weak solution of (\ref{3.1}), we
apply (\ref{3.4}) to $u_{t}$ and derive
\begin{equation}
\sum\limits_{\left\vert I\right\vert \leq k-2}\iint\nolimits_{Q_{R/2^{k-1}}%
}{{\left\vert {X_{I}u_{t}}\right\vert ^{2}dz}}\leq\frac{c}{{R^{2(k-2)}}}%
\iint\nolimits_{Q_{R/2}}{\left\vert {u_{t}}\right\vert ^{2}dz}\leq\frac
{c}{{R^{2k}}}\iint\nolimits_{Q_{R}}{\left\vert u\right\vert ^{2}dz}.
\label{3.14}%
\end{equation}
Inserting (\ref{3.13}) and (\ref{3.14}) into (\ref{3.12}), it shows from
(\ref{3.4}) that%
\begin{align}
II  &  \leq\frac{c}{{R^{2k}}}\iint\nolimits_{Q_{R/2^{k-1}}}{\left\vert
u\right\vert ^{2}dz}+\frac{c}{{R^{2k}}}\iint\nolimits_{Q_{R}}{\left\vert
u\right\vert ^{2}dz} \leq\frac{c}{{R^{2k}}}\iint\nolimits_{Q_{R}}{\left\vert
u\right\vert ^{2}dz}. \label{3.15}%
\end{align}

Now let us handle $I$. Since $c_{\alpha}\in C_{0}^{\infty}$, $\left\vert
{X_{I}c_{\alpha}}\right\vert $ is bounded. By properties of $\xi\left(
x\right)  $ and $\eta\left(  t\right)  $,%
\begin{align}
I  &  =c\int_{I_{R/2^{k-1}}}{(\sum\limits_{\left\vert I\right\vert \leq
k-2}{\int_{B_{R/2^{k-1}}}{\left\vert {a_{ij}^{\alpha\beta}\eta X_{I}%
(c_{\alpha}X_{\beta}(\xi u)^{j})}\right\vert ^{2}}dx)}}dt\nonumber\\
&  \leq c\sum\limits_{\left\vert I\right\vert \leq k-2}\iint
\nolimits_{Q_{R/2^{k-1}}}{{[\left\vert {\xi X_{\beta}u^{j}X_{I}c_{\alpha}%
}\right\vert ^{2}+\left\vert {c_{\alpha}X_{\beta}u^{j}X_{I}\xi}\right\vert
^{2}+\left\vert {c_{\alpha}\xi X_{I}X_{\beta}u^{j}}\right\vert ^{2}]dz}%
}\nonumber\\
&  +c\sum\limits_{\left\vert I\right\vert \leq k-2}\iint
\nolimits_{Q_{R/2^{k-1}}}{{[\left\vert {u^{j}X_{\beta}\xi X_{I}c_{\alpha}%
}\right\vert ^{2}+\left\vert {c_{\alpha}X_{\beta}\xi X_{I}u^{j}}\right\vert
^{2}+\left\vert {c_{\alpha}u^{j}X_{I}X_{\beta}\xi}\right\vert ^{2}]dz}%
}\nonumber\\
&  \leq c\sum\limits_{\left\vert I\right\vert \leq k-2}\iint
\nolimits_{Q_{R/2^{k-1}}}{{[\left\vert {X_{\beta}u^{j}}\right\vert ^{2}%
+\frac{c}{{R^{2\left\vert I\right\vert }}}\left\vert {X_{\beta}u^{j}%
}\right\vert ^{2}+\left\vert {X_{I}X_{\beta}u^{j}}\right\vert ^{2}]dz}%
}\nonumber\\
&  +c\sum\limits_{\left\vert I\right\vert \leq k-2}\iint
\nolimits_{Q_{R/2^{k-1}}}{{[\frac{1}{{R^{2}}}\left\vert {u^{j}}\right\vert
^{2}+\frac{1}{{R^{2}}}\left\vert {X_{I}u^{j}}\right\vert ^{2}+\frac
{c}{{R^{2(\left\vert I\right\vert +1)}}}\left\vert {u^{j}}\right\vert ^{2}%
]dz}}. \label{3.16}%
\end{align}
By the assertion for $\left\vert I\right\vert \leq k-1$,
\[
\sum\limits_{\left\vert I\right\vert \leq k-2}\iint\nolimits_{Q_{R/2^{k-1}}%
}{{\left\vert {X_{I}X_{\beta}u^{j}}\right\vert ^{2}dz}}\leq\frac
{c}{{R^{2(k-1)}}}\iint\nolimits_{Q_{R}}{\left\vert u\right\vert ^{2}dz}%
\leq\frac{c}{{R^{2k}}}\iint\nolimits_{Q_{R}}{\left\vert u\right\vert ^{2}dz}%
\]
and%
\[
\sum\limits_{\left\vert I\right\vert \leq k-2}\iint\nolimits_{Q_{R/2^{k-1}}%
}{{\left\vert {X_{I}u^{j}}\right\vert ^{2}dz}}\leq\frac{c}{{R^{2(k-2)}}}%
\iint\nolimits_{Q_{R/2}}{\left\vert u\right\vert ^{2}dz}\leq\frac
{c}{{R^{2(k-2)}}}\iint\nolimits_{Q_{R}}{\left\vert u\right\vert ^{2}dz}.
\]
Inserting the above two inequalities into (\ref{3.16}) yields%
\begin{align}
I  &  \leq\frac{c}{{R^{2}}}\iint\nolimits_{Q_{R/2^{k-2}}}{\left\vert
u\right\vert ^{2}dz}+\frac{c}{{R^{2(k-1)}}}\iint\nolimits_{Q_{R}}{\left\vert
u\right\vert ^{2}dz}+\frac{c}{{R^{2k}}}\iint\nolimits_{Q_{R}}{\left\vert
u\right\vert ^{2}dz}\nonumber\\
&  \leq\frac{c}{{R^{2k}}}\iint\nolimits_{Q_{R}}{\left\vert u\right\vert
^{2}dz}. \label{3.17}%
\end{align}
Putting (\ref{3.15}) and (\ref{3.17}) into (\ref{3.11}), we get
\[
\sum\limits_{\left\vert I\right\vert =k}\iint\nolimits_{Q_{R/2^{k}}%
}{{\left\vert {X_{I}u}\right\vert ^{2}dz}}\leq\frac{c}{{R^{2k}}}%
\iint\nolimits_{Q_{R}}{\left\vert u\right\vert ^{2}dz},
\]
hence (\ref{3.9}) is proved.

The proof of (\ref{3.10}) is easy. In fact, since $\partial_{t}^{m}u$ is also
a weak solution of (\ref{3.1}), it shows by applying (\ref{3.9}) to
$\partial_{t}^{m}u$ and noting (\ref{3.4}) that%
\begin{align*}
&  \sum\limits_{\left\vert I\right\vert +2m\leq k}\iint\nolimits_{Q_{R/2^{k}}%
}{{\left\vert {X_{I}\partial_{t}^{m}u}\right\vert ^{2}dxdt}}\\
&  \leq\sum\limits_{\left\vert I\right\vert +2m\leq k}{\frac{c}%
{{R^{2\left\vert I\right\vert }}}\iint\nolimits_{Q_{R/2^{2m}}}{\left\vert
{\partial_{t}^{m}u}\right\vert ^{2}dxdt}}\\
&  \leq\sum\limits_{\left\vert I\right\vert +2m\leq k}{\frac{c}%
{{R^{2\left\vert I\right\vert }}}\frac{c}{{R^{4m}}}\iint\nolimits_{Q_{R}%
}{\left\vert u\right\vert ^{2}dxdt}}\\
&  \leq\frac{c}{{R^{2k}}}\iint\nolimits_{Q_{R}}{\left\vert u\right\vert
^{2}dxdt}.
\end{align*}

\end{proof}

\begin{lemma}
\textbf{\label{l3.4}} (Sobolev-Poincar\'{e} inequality, see \cite{b7},
\cite{b12} and references therein). For any open set $\Omega^{\prime}$,
$\bar{\Omega}^{\prime}\Subset\Omega$, there exist positive constants $R_{0}$
and $c$ such that for any $x_{0}\in\Omega^{\prime}$, $0<R\leq R_{0}$, $u\in
C^{\infty}(\bar{B}_{R})$,
\[
(\frac{1}{{\left\vert {B_{R}}\right\vert }}\int_{B_{R}}{\left\vert {u-u_{R}%
}\right\vert ^{q^{\prime}}}dx)^{\frac{1}{{q^{\prime}}}}\leq cR(\frac
{1}{{\left\vert {B_{R}}\right\vert }}\int_{B_{R}}{\sum\limits_{i=1}%
^{q}{\left\vert {X_{i}u}\right\vert ^{p^{\prime}}}}dx)^{\frac{1}{{p^{\prime}}%
}},
\]
where $1<p^{\prime}<Q$, $1\leq q^{\prime}<\frac{{p^{\prime}Q}}{{Q-p^{\prime}}%
}$, $u_{R}\left(  t\right)  =\frac{1}{{\left\vert {B_{R}}\right\vert }}%
\int_{B_{R}}{u\left(  {x,t}\right)  }dx$, $R_{0}$ and $c$ depend on
$\Omega^{\prime}$ and $\Omega$.

If $u\in C_{0}^{\infty}(B_{R})$, then for all $1\leq q^{\prime}\leq
\frac{{p^{\prime}Q}}{{Q-p^{\prime}}}$,
\begin{equation}
(\frac{1}{{\left\vert {B_{R}}\right\vert }}\int_{B_{R}}{\left\vert
u\right\vert ^{q^{\prime}}}dx)^{\frac{1}{{q^{\prime}}}}\leq cR(\frac
{1}{{\left\vert {B_{R}}\right\vert }}\int_{B_{R}}{\sum\limits_{i=1}%
^{q}{\left\vert {X_{i}u}\right\vert ^{p^{\prime}}}}dx)^{\frac{1}{{p^{\prime}}%
}}. \label{3.18}%
\end{equation}
In particular, if $p^{\prime}=q^{\prime}=2$, then
\begin{equation}
\int_{B_{R}}{\left\vert {u}\right\vert ^{2}}dx \leq c {R^{2}} \int_{B_{R}%
}{\sum\limits_{i=1}^{q}{\left\vert {X_{i}u}\right\vert ^{2}}}dx; \label{3.19}%
\end{equation}
if $p^{\prime}=2$, $q^{\prime}=\frac{{2Q}}{{Q-2}}$, then
\begin{equation}
(\int_{B_{R}}{\left\vert u\right\vert ^{\frac{{2Q}}{{Q-2}}}}dx)^{\frac{{Q-2}%
}{{2Q}}}\leq c(\int_{B_{R}}{\left\vert {Xu}\right\vert ^{2}}dx)^{\frac{1}{2}}.
\label{3.20}%
\end{equation}

\end{lemma}

\begin{lemma}
\label{l3.5} Let $u\in V_{2}(Q_{T})$ be a weak solution of (\ref{3.1}) in
$Q_{T}$ and $Q_{R}\subset Q_{T}$. Then for any $0\leq\rho\leq R$,
\[
\iint\nolimits_{Q_{\rho}}{\left\vert u\right\vert ^{2}dz}\leq C(\frac{\rho}%
{R})^{Q+2}\iint\nolimits_{Q_{R}}{\left\vert u\right\vert ^{2}dz}.
\]

\end{lemma}

\begin{proof}
Let $k_{1}$ and $k_{2}$ be fixed integers such that $k_{1}>\frac{Q}{2}$ and
$k_{2}>1$. If $\rho\geq\frac{R}{{2^{k_{1}+k_{2}+2}}}$, then the conclusion is
obvious. If $\rho<\frac{R}{{2^{k_{1}+k_{2}+2}}}$, then by (\ref{2.1}),
(\ref{3.7}) and (\ref{3.8}),%
\begin{align*}
\iint\nolimits_{Q_{\rho}}{\left\vert u\right\vert ^{2}dxdt}  &  \leq
\int_{I_{\rho}}{\left\vert {B_{\rho}}\right\vert
\mathop {\sup }\limits_{B_{R/2^{k_{1}+k_{2}+2}}}\left\vert u\right\vert ^{2}%
}dt\\
&  \leq c\left\vert {B_{\rho}}\right\vert \int_{I_{\rho}}{(\sum
\limits_{\left\vert I\right\vert \leq k_{1}}{\left\vert {B_{R/2^{k_{1}+k_{2}}%
}}\right\vert ^{-1}R^{2\left\vert I\right\vert }\int_{B_{R/2^{k_{1}+k_{2}}}%
}{\left\vert {X_{I}u}\right\vert ^{2}}dx)}}dt\\
&  \leq c\frac{{\left\vert {B_{\rho}}\right\vert }}{{\left\vert {B_{R}%
}\right\vert }}\sum\limits_{\left\vert I\right\vert \leq k_{1}}{R^{2\left\vert
I\right\vert }\int_{B_{R/2^{k_{1}+k_{2}}}}{\left\vert {I_{\rho}}\right\vert
\mathop {\sup }\limits_{I_{\rho}}\left\vert {X_{I}u}\right\vert ^{2}}dx}\\
&  \leq c\rho^{2}\frac{{\left\vert {B_{\rho}}\right\vert }}{{\left\vert
{B_{R}}\right\vert }}\sum\limits_{\left\vert I\right\vert \leq k_{1}%
}{R^{2\left\vert I\right\vert }\int_{B_{R/2^{k_{1}+k_{2}}}}{\sum
\limits_{2m\leq k_{2}}{R^{4m-2}\int_{I_{R/2^{k_{1}+k_{2}}}}{\left\vert
{\partial_{t}^{m}X_{I}u}\right\vert ^{2}}dtdx}}}\\
&  \leq c\left(  {\frac{\rho}{R}}\right)  ^{2}\frac{{\left\vert {B_{\rho}%
}\right\vert }}{{\left\vert {B_{R}}\right\vert }}\sum\limits_{\left\vert
I\right\vert +2m\leq k_{1}+k_{2}}{R^{2(\left\vert I\right\vert +2m)}%
\iint\nolimits_{Q_{R/2^{k_{1}+k_{2}}}}{\left\vert {\partial_{t}^{m}X_{I}%
u}\right\vert ^{2}dz}}.
\end{align*}
Applying (\ref{3.10}) leads to
\[
\iint\nolimits_{Q_{\rho}}{\left\vert u\right\vert ^{2}dxdt}\leq c\left(
{\frac{\rho}{R}}\right)  ^{2}\frac{{\left\vert {B_{\rho}}\right\vert }%
}{{\left\vert {B_{R}}\right\vert }}\iint\nolimits_{Q_{R}}{\left\vert
u\right\vert ^{2}dz}\leq c\left(  {\frac{\rho}{R}}\right)  ^{Q+2}%
\iint\nolimits_{Q_{R}}{\left\vert u\right\vert ^{2}dz},
\]
where we have used the definition of $Q$ and the fact that $\left\vert {B_{R}%
}\right\vert $ can be approximated by some polynomial in $R$, see
\cite{b7},\cite{b14}.
\end{proof}

\begin{lemma}
\label{l3.6} Suppose that $u\in V_{2}(Q_{T})$ is a weak solution of
(\ref{3.1}), $Q_{R}(z_{0})\subset Q_{T}$ and $u=0$ on $\partial_{p}Q_{R}$.
Then for any $0\leq\rho\leq R$, it follows
\[
\iint\nolimits_{Q_{\rho}}{\left\vert {Xu}\right\vert ^{2}dxdt}\leq c\left(
{\frac{\rho}{R}}\right)  ^{Q+2}\iint\nolimits_{Q_{R}}{\left\vert
{Xu}\right\vert ^{2}dxdt}.
\]

\end{lemma}

\begin{proof}
Let $k_{1}$ and $k_{2}$ be fixed integers such that $k_{1}>\frac{Q}{2}$ and
$k_{2}>1$. If $\rho\geq\frac{R}{{2^{k_{1}+k_{2}+3}}}$, then the conclusion
holds; if $\rho<\frac{R}{{2^{k_{1}+k_{2}+3}}}$, then by (\ref{3.7}) and
(\ref{3.8}),%
\begin{align*}
\iint\nolimits_{Q_{\rho}}{\left\vert {X_{i}u}\right\vert ^{2}dxdt} &  \leq
\int_{I_{\rho}}{\left\vert {B_{\rho}}\right\vert
\mathop {\sup }\limits_{B_{R/2^{k_{1}+k_{2}+3}}}\left\vert {X_{i}u}\right\vert
^{2}}dt\\
&  \leq c\left\vert {B_{\rho}}\right\vert \int_{I_{\rho}}{(\sum
\limits_{\left\vert I\right\vert \leq k_{1}}{\left\vert {B_{R/2^{k_{1}%
+k_{2}+1}}}\right\vert ^{-1}R^{2\left\vert I\right\vert }\int_{B_{R/2^{k_{1}%
+k_{2}+1}}}{\left\vert {X_{I}X_{i}u}\right\vert ^{2}}dx)}}dt\\
&  \leq c\frac{{\left\vert {B_{\rho}}\right\vert }}{{\left\vert {B_{R}%
}\right\vert }}\sum\limits_{\left\vert I\right\vert \leq k_{1}}{R^{2\left\vert
I\right\vert }\int_{B_{R/2^{k_{1}+k_{2}+1}}}{\left\vert {I_{\rho}}\right\vert
\mathop {\sup }\limits_{I_{\rho}}\left\vert {X_{I}X_{i}u}\right\vert ^{2}}%
dx}\\
&  \leq c\rho^{2}\frac{{\left\vert {B_{\rho}}\right\vert }}{{\left\vert
{B_{R}}\right\vert }}\sum\limits_{\left\vert I\right\vert \leq k_{1}%
}{R^{2\left\vert I\right\vert }\int_{B_{R/2^{k_{1}+k_{2}+1}}}{\sum
\limits_{2m\leq k_{2}}{R^{4m-2}\int_{I_{R/2^{k_{1}+k_{2}+1}}}{\left\vert
{\partial_{t}^{m}X_{I}X_{i}u}\right\vert ^{2}}dt}}}dx\\
&  \leq c\left(  {\frac{\rho}{R}}\right)  ^{2}\frac{{\left\vert {B_{\rho}%
}\right\vert }}{{\left\vert {B_{R}}\right\vert }}\sum\limits_{\left\vert
I\right\vert +2m\leq k_{1}+k_{2}}{R^{2(\left\vert I\right\vert +2m)}%
\iint\nolimits_{Q_{R/2^{k_{1}+k_{2}+1}}}{\left\vert {\partial_{t}^{m}%
X_{I}X_{i}u}\right\vert ^{2}dz}}.
\end{align*}
In virtue of (\ref{3.10}) and (\ref{3.18}),%
\begin{align*}
\iint\nolimits_{Q_{\rho}}{\left\vert {X_{i}u}\right\vert ^{2}dxdt} &  \leq
c\left(  {\frac{\rho}{R}}\right)  ^{2}\frac{{\left\vert {B_{\rho}}\right\vert
}}{{\left\vert {B_{R}}\right\vert }}\sum\limits_{\left\vert I\right\vert
+2m\leq k}{R^{2(\left\vert I\right\vert +2m)}R^{-2(\left\vert I\right\vert
+2m+1)}\iint\nolimits_{Q_{R}}{\left\vert u\right\vert ^{2}dz}}\\
&  \leq c\left(  {\frac{\rho}{R}}\right)  ^{2}\frac{{\left\vert {B_{\rho}%
}\right\vert }}{{\left\vert {B_{R}}\right\vert }}\frac{1}{{R^{2}}}%
\iint\nolimits_{Q_{R}}{\left\vert u\right\vert ^{2}dz},\\
&  \leq c\left(  {\frac{\rho}{R}}\right)  ^{2}\frac{{\left\vert {B_{\rho}%
}\right\vert }}{{\left\vert {B_{R}}\right\vert }}\iint\nolimits_{Q_{R}%
}{\left\vert {Xu}\right\vert ^{2}dz}.
\end{align*}
Similarly to Lemma 3.6, we end the proof.
\end{proof}

We need to define a parabolic distance $d_{p}$ corresponding to $d_{X}$. For
$(x,t),(y,s)\in Q_{T}$, set
\[
d_{p}\left(  (x,t),(y,s)\right)  =({d_{X}(x,y)^{2}+\left\vert
{t-s}\right\vert })^{\frac{1}{2}}.
\]
Denote a ball with respect to the distance $d_{p}$ by
\[
B_{p}(\left(  x_{0},t_{0}\right)  ,R)=\{\left(  x,t\right)  \in Q_{T}%
:d_{p}(\left(  x_{0},t_{0}\right)  ,\left(  x,t\right)  )<R\}.
\]

An important fact is that $B_{p}(\left(  x_{0},t_{0}\right)  ,R)$ is
a homogeneous space (see [8], [1, Proposition 3.8]).
According to it and
\[
Q_{R}\left(  z\right)  \subset B_{p}\left(  z,2R\right)  \subset Q_{2R}\left(
z\right)  ,
\]
we immediately know that the reverse H\"{o}lder inequality in
\cite{b8} (or \cite{b21}) is true for parabolic cylinders.

\begin{lemma}
\label{l3.7}  Let
$g\geq0$ on $Q_{T}$. If for some $\hat{q}>1$ such that for any $Q_{4R}\subset
Q_{T}$,
\begin{equation}
\frac{1}{{\left\vert {Q_{R}}\right\vert }}\iint\nolimits_{Q_{R}}{g^{\hat{q}%
}dxdt}\leq b\left(  {\frac{1}{{\left\vert {Q_{4R}}\right\vert }}%
\iint\nolimits_{Q_{4R}}{gdxdt}}\right)  ^{\hat{q}}+\theta{\iint
\nolimits_{Q_{4R}}g^{\hat{q}}dxdt}.\label{3.21}%
\end{equation}
Then there exist positive constants $b>1$ and $\theta_{0}=\theta_{0}(\hat
{q},Q_{T})$ such that if $\theta<\theta_{0}$, then $g\in L_{loc}%
^{p}\left(  {Q_{T}}\right)  $ for any $p\in\left(  {\hat{q},\hat
{q}+\varepsilon}\right]  $. Moreover, it holds
\[
\left(  {\frac{1}{{\left\vert {Q_{R}}\right\vert }}\iint\nolimits_{Q_{R}%
}{g^{p}dxdt}}\right)  ^{\frac{1}{p}}\leq c\left(  {\frac{1}{{\left\vert
{Q_{4R}}\right\vert }}\iint\nolimits_{Q_{4R}}{g^{\hat{q}}dxdt}}\right)
^{\frac{1}{{\hat{q}}}},
\]
where the positive constants $c$ and $\varepsilon$ depend only on $b,\hat
{q},\theta$ and $Q$.
\end{lemma}

\begin{theorem}
\label{T3.1} Let $Q_{R}\subset Q_{4R}\subset Q_{T}$ and $u\in V_{2}(Q_{T})$ be
a weak solution of (\ref{3.1}) in $Q_{T}$ and $u=0$ on $\partial_{p}Q_{4R}$.
Then there exists a constant $s>2$ such that $Xu\in L_{loc}^{s}(Q_{T})$.
Moreover, the following inequality holds
\[
(\frac{1}{{\left\vert {Q_{R}}\right\vert }}\iint\nolimits_{Q_{R}}{\left\vert
{Xu}\right\vert }^{s}{dz})^{1/s}\leq C(\frac{1}{{\left\vert {Q_{4R}%
}\right\vert }}{\iint\nolimits_{Q_{4R}}\left\vert {Xu}\right\vert ^{2}%
dz})^{1/2}.
\]

\end{theorem}

\begin{proof}
Set $2^{\ast}=\frac{{2Q}}{{Q-2}}$ and $\tilde{q}=\frac{{2Q}}{{Q+2}}$. Note%
\begin{align}
&  \iint\nolimits_{Q_{2R}}{\left\vert {u\left(  t\right)  }\right\vert ^{2}%
dz}\nonumber\\
&  \leq\mathop {\sup }\limits_{I_{2R}}\left(  {\int_{B_{2R}}{\left\vert
{u\left(  t\right)  }\right\vert ^{2}dx}}\right)  ^{\frac{1}{2}}\cdot\left(
{\int_{I_{2R}}{\left(  {\int_{B_{2R}}{\left\vert {u\left(  t\right)
}\right\vert ^{2}dx}}\right)  }^{\frac{1}{2}}dt}\right)  ,\label{3.22}%
\end{align}
and denote
\[
A\equiv\mathop {\sup }\limits_{I_{2R}}\left(  {\int_{B_{2R}}{\left\vert
{u\left(  t\right)  }\right\vert ^{2}dx}}\right)  ^{\frac{1}{2}},B\equiv
\int_{I_{2R}}{\left(  {\int_{B_{2R}}{\left\vert {u\left(  t\right)
}\right\vert ^{2}dx}}\right)  }^{\frac{1}{2}}dt.
\]
Now we estimate $A$ and $B$, respectively. By (\ref{3.2}) and (\ref{3.19}),
\begin{equation}
A\leq\frac{c}{R}(\iint\nolimits_{Q_{4R}}{\left\vert {u}\right\vert ^{2}%
dz})^{\frac{1}{2}}\leq c(\iint\nolimits_{Q_{4R}}{\left\vert {Xu}\right\vert
^{2}dz})^{\frac{1}{2}}.\label{3.23}%
\end{equation}
To B, we have by (\ref{3.18}) and (\ref{3.20}) that
\begin{align}
B &  \leq\int_{I_{2R}}{\left(  {\int_{B_{2R}}{\left\vert {u\left(  t\right)
}\right\vert ^{\tilde{q}}dx}}\right)  }^{\frac{1}{{2\tilde{q}}}}\left(
{\int_{B_{2R}}{\left\vert {u\left(  t\right)  }\right\vert ^{2^{\ast}}dx}%
}\right)  ^{\frac{1}{{2\cdot2^{\ast}}}}dt\nonumber\\
&  \leq c\int_{I_{2R}}{\left(  R^{\tilde{q}}\int_{B_{2R}}{\left\vert
{Xu}\right\vert ^{\tilde{q}}dx}\right)  }^{\frac{1}{{2\tilde{q}}}}\left(
{\int_{B_{2R}}{\left\vert {u}\right\vert ^{2^{\ast}}dx}}\right)  ^{\frac
{1}{{2\cdot2^{\ast}}}}dt\nonumber\\
&  \leq cR^{\frac{1}{2}}\int_{I_{2R}}{\left(  {\int_{B_{2R}}{\left\vert
{Xu}\right\vert ^{\tilde{q}}dx}}\right)  }^{\frac{1}{{2\tilde{q}}}}\left(
{\int_{B_{2R}}{\left\vert {Xu}\right\vert ^{2}dx}}\right)  ^{\frac{1}{4}%
}dt\nonumber\\
&  \leq cR^{\frac{1}{2}}\left(  \iint\nolimits_{Q_{2R}}{{\left\vert
{Xu}\right\vert ^{\tilde{q}}dz}}\right)  ^{\frac{1}{{2\tilde{q}}}}\left(
{\int_{I_{2R}}{\left(  {\int_{B_{2R}}{\left\vert {Xu}\right\vert ^{2}dx}%
}\right)  }^{\frac{1}{2}\frac{{\tilde{q}}}{{2\tilde{q}-1}}}dt}\right)
^{\frac{{2\tilde{q}-1}}{{2\tilde{q}}}}\nonumber\\
&  \leq cR^{\frac{3}{2}-\frac{1}{Q}}\left(  \iint\nolimits_{Q_{4R}%
}{{\left\vert {Xu}\right\vert ^{\tilde{q}}dz}}\right)  ^{\frac{1}{{2\tilde{q}%
}}}\left(  \iint\nolimits_{Q_{4R}}{{\left\vert {Xu}\right\vert ^{2}dz}%
}\right)  ^{\frac{1}{4}}.\label{3.24}%
\end{align}
Inserting (\ref{3.23}) and (\ref{3.24}) into (\ref{3.22}) and using Young's
inequality,%
\begin{align*}
&  \iint\nolimits_{Q_{2R}}{\left\vert {u\left(  t\right)  }\right\vert ^{2}%
dz}\\
&  \leq cR^{\frac{3}{2}-\frac{1}{Q}}\left(  \iint\nolimits_{Q_{4R}%
}{{\left\vert {Xu}\right\vert ^{2}dz}}\right)  ^{\frac{3}{4}}\cdot\left(
\iint\nolimits_{Q_{4R}}{{\left\vert {Xu}\right\vert ^{\tilde{q}}dz}}\right)
^{\frac{1}{{2\tilde{q}}}}\\
&  \leq\varepsilon R^{2}\iint\nolimits_{Q_{4R}}{\left\vert {Xu}\right\vert
^{2}dz}+C(\varepsilon)R^{-\frac{4}{Q}}\left(  \iint\nolimits_{Q_{4R}%
}{{\left\vert {Xu}\right\vert ^{\tilde{q}}dz}}\right)  ^{\frac{2}{{\tilde{q}}%
}}.
\end{align*}
Returning to (\ref{3.2}) and using the above inequality lead to%
\begin{align*}
&  \frac{1}{{\left\vert {Q_{R}}\right\vert }}\iint\nolimits_{Q_{R}}{\left\vert
{Xu}\right\vert ^{2}dz}\\
&  \leq\frac{c}{{R^{2}}}\frac{1}{{\left\vert {Q_{R}}\right\vert }}%
\iint\nolimits_{Q_{2R}}{\left\vert {u\left(  t\right)  }\right\vert ^{2}dz}\\
&  \leq\varepsilon\frac{1}{{\left\vert {Q_{4R}}\right\vert }}\iint
\nolimits_{Q_{4R}}{\left\vert {Xu}\right\vert ^{2}dz}+C(\varepsilon
)\frac{{\left\vert {Q_{4R}}\right\vert ^{2/\tilde{q}}}}{{\left\vert {Q_{R}%
}\right\vert }}R^{-\frac{4}{Q}-2}\left(  {\frac{1}{{\left\vert {Q_{4R}%
}\right\vert }}\iint\nolimits_{Q_{4R}}{\left\vert {Xu}\right\vert ^{\tilde{q}%
}dz}}\right)  ^{\frac{2}{{\tilde{q}}}}\\
&  \leq\varepsilon\frac{1}{{\left\vert {Q_{4R}}\right\vert }}\iint
\nolimits_{Q_{4R}}{\left\vert {Xu}\right\vert ^{2}dz}+C(\varepsilon)\left(
{\frac{1}{{\left\vert {Q_{4R}}\right\vert }}\iint\nolimits_{Q_{4R}}{\left\vert
{Xu}\right\vert ^{\tilde{q}}dz}}\right)  ^{\frac{2}{{\tilde{q}}}}.
\end{align*}
Let $g=\left\vert {Xu}\right\vert ^{\tilde{q}}$, $\hat{q}=\frac{2}{{\tilde{q}%
}}=\frac{{Q+2}}{Q}>1$ , $\theta=\varepsilon$. The previous inequality is of the
form
\[
\frac{1}{{\left\vert {Q_{R}}\right\vert }}\iint\nolimits_{Q_{R}}{g^{\hat{q}%
}dz}\leq\theta\frac{1}{{\left\vert {Q_{4R}}\right\vert }}\iint
\nolimits_{Q_{4R}}{g^{\hat{q}}dz}+C(\varepsilon)\left(  {\frac{1}{{\left\vert
{Q_{4R}}\right\vert }}\iint\nolimits_{Q_{4R}}{gdz}}\right)  ^{\hat{q}}.
\]
Due to Lemma \ref{l3.7}, there exists $\varepsilon>0$ such that
for any $p\in\lbrack\hat{q},\hat{q}+\varepsilon)$,
\[
\left(  {\frac{1}{{\left\vert {Q_{R}}\right\vert }}\iint\nolimits_{Q_{R}%
}{\left\vert {Xu}\right\vert ^{p\tilde{q}}dz}}\right)  ^{\frac{1}{p}}\leq
c\left(  {\frac{1}{{\left\vert {Q_{4R}}\right\vert }}\iint\nolimits_{Q_{4R}%
}{\left\vert {Xu}\right\vert ^{2}dz}}\right)  ^{\frac{{\tilde{q}}}{2}},
\]
Denoting $s=p\tilde{q}\in\lbrack2,2+\varepsilon)$, the proof is finished.
\end{proof}

\begin{remark}
\label{R3.2}It is not hard to find that the conclusion of Theorem \ref{T3.1}
is still true for the homogeneous parabolic system with variable coefficients,
when we check carefully the above proof. It will be useful in Section 4.
\end{remark}

\section{Proof of Theorem 1.1}

In this section we will prove Theorem \ref{T1.1}. First step is
to establish the following.

\begin{theorem}
\label{T4.1} Let $u\in V_{2}(Q_{T})$ be a weak solution of
\[
u_{t}^{i}+X_{\alpha}^{\ast}(a_{ij}^{\alpha\beta}(z)X_{\beta}u^{j})=0,
\]
in $Q_{T}$. Suppose coefficients $a_{ij}^{\alpha\beta}(z)\in VMO(Q_{T})$ and
satisfy (\ref{1.2}). Then for any $0<\mu<Q+2$, there exist positive constants
$R_{0}$ and $c$ such that for any $\rho\leq R\leq\frac{1}{2}\min
(R_{0},dist(z_{0},\partial_{p}Q_{T}))$, it holds
\[
{\iint\nolimits_{Q_{\rho}}\left\vert {Xu}\right\vert ^{2}dxdt}\leq c\left(
{\frac{\rho}{R}}\right)  ^{\mu}{\iint\nolimits_{Q_{R}}\left\vert
{Xu}\right\vert ^{2}dxdt},
\]
where $R_{0}$ and $c$ depend on $Q$, $\mu$, $\Lambda$ and the $VMO$ modulus of
$a_{ij}^{\alpha\beta}$.
\end{theorem}

\begin{proof}
Let $w$ be a weak solution of the following system
\begin{equation}
\left\{
\begin{array}
[c]{l}%
w_{t}^{i}+X_{\alpha}^{\ast}((a_{ij}^{\alpha\beta})_{z_{0},R})X_{\beta}%
w^{j})=0,\mathrm{{in}}Q_{R},\\
\text{ \ \ \ \ \ }w=u,\text{ \ \ \ \ \ \ \ \ \ }\mathrm{{on}}\text{ }%
\partial_{p}Q_{R},
\end{array}
\right.  \label{4.1}%
\end{equation}
where $z_{0}$ is a fixed point in $Q_{R}$, $(a_{ij}^{\alpha\beta})_{z_{0}%
,R}=\frac{1}{{\left\vert {Q_{T}\cap Q_{R}}\right\vert }}{\iint\nolimits_{Q_{T}%
\cap Q_{R}}a_{ij}^{\alpha\beta}(z)dz}$. Then $v=u-w$ satisfies
\begin{equation}
\left\{
\begin{array}
[c]{l}%
v_{t}^{i}+X_{\alpha}^{\ast}((a_{ij}^{\alpha\beta})_{z_{0},R}X_{\beta}%
v^{j})=X_{\alpha}^{\ast}(((a_{ij}^{\alpha\beta})_{z_{0},R}-a_{ij}^{\alpha
\beta}(z))X_{\beta}u^{j}),\mathrm{{in}}\text{ }Q_{R},\\
\text{ \ \ \ \ \ \ \ \ \ \ \ \ \ \ \ \ \ \ \ }v=0,\text{
\ \ \ \ \ \ \ \ \ \ \ \ \ \ \ \ \ }\mathrm{{on}}\text{ }\partial_{p}Q_{R}.
\end{array}
\right.  \label{4.2}%
\end{equation}
Multiplying both sides of (\ref{4.2}) by $v^{i}$ and integrating by parts on
$Q_{R}$,
\begin{align*}
&  {\iint\nolimits_{Q_{R}}(v_{t}^{i}v^{i}+(a_{ij}^{\alpha\beta})_{z_{0}%
,R}X_{\beta}v^{j}X_{\alpha}v^{i})dz}\\
&  ={\iint\nolimits_{Q_{R}}(((a_{ij}^{\alpha\beta})_{z_{0},R}-a_{ij}%
^{\alpha\beta}(z))X_{\beta}u^{j}X_{\alpha}v^{i})dz}\\
&  \leq{\iint\nolimits_{Q_{R}}\left\vert {a_{ij}^{\alpha\beta}(z)-(a_{ij}%
^{\alpha\beta})_{z_{0},R}}\right\vert \left\vert {Xu}\right\vert \left\vert
{Xv}\right\vert dz}\\
&  \leq C_{\varepsilon}{\iint\nolimits_{Q_{R}}\left\vert {a_{ij}^{\alpha\beta
}(z)-(a_{ij}^{\alpha\beta})_{z_{0},R}}\right\vert ^{2}\left\vert
{Xu}\right\vert ^{2}dz}+\varepsilon{\iint\nolimits_{Q_{R}}\left\vert
{Xv}\right\vert ^{2}dz}.
\end{align*}
Noting ${\iint\nolimits_{Q_{R}}v_{t}^{i}v^{i}dz}=\int_{B_{R}}{dx\int
_{t_{0}-R^{2}}^{t_{0}}{v^{i}dv^{i}}}\geq0$ and (\ref{1.2}), it follows
\begin{equation}
{\iint\nolimits_{Q_{R}}\left\vert {Xv}\right\vert ^{2}dz}\leq C_{\varepsilon
}{\iint\nolimits_{Q_{R}}\left\vert {a_{ij}^{\alpha\beta}(z)-(a_{ij}%
^{\alpha\beta})_{z_{0},R}}\right\vert ^{2}\left\vert {Xu}\right\vert ^{2}%
dz}.\label{4.3}%
\end{equation}
From $a_{ij}^{\alpha\beta}\in VMO$, we see that for any $\varepsilon>0$, there
exists $R_{0}>0$ such that for any $R\leq R_{0}$,
\[
\left(  {\frac{1}{{\left\vert {Q_{R}}\right\vert }}\iint\nolimits_{Q_{R}%
}{\left\vert {a_{ij}^{\alpha\beta}(z)-(a_{ij}^{\alpha\beta})_{z_{0},R}%
}\right\vert ^{\frac{{2s}}{{s-2}}}dz}}\right)  ^{\frac{{s-2}}{s}}<\varepsilon
\]
and%
\begin{align*}
&  {\iint\nolimits_{Q_{R}}\left\vert {a_{ij}^{\alpha\beta}(z)-(a_{ij}%
^{\alpha\beta})_{z_{0},R}}\right\vert ^{2}\left\vert {Xu}\right\vert ^{2}dz}\\
&  \leq\left\vert {Q_{R}}\right\vert (\frac{1}{{\left\vert {Q_{R}}\right\vert
}}{\iint\nolimits_{Q_{R}}\left\vert {a_{ij}^{\alpha\beta}(z)-(a_{ij}%
^{\alpha\beta})_{z_{0},R}}\right\vert ^{\frac{{2s}}{{s-2}}}dz})^{\frac{{s-2}%
}{s}}(\frac{1}{{\left\vert {Q_{R}}\right\vert }}{\iint\nolimits_{Q_{R}%
}\left\vert {Xu}\right\vert ^{s}dz})^{2/s}\\
&  \leq\varepsilon\left\vert {Q_{R}}\right\vert (\frac{1}{{\left\vert {Q_{R}%
}\right\vert }}{\iint\nolimits_{Q_{R}}\left\vert {Xu}\right\vert ^{s}%
dz})^{2/s}\\
&  \leq\varepsilon{\iint\nolimits_{Q_{4R}}\left\vert {Xu}\right\vert ^{2}dz},
\end{align*}
where we have used Theorem \ref{T3.1}, Remarks \ref{R3.1} and \ref{R3.2}.

Inserting the above inequality into (\ref{4.3}), we immediately get
\[
{\iint\nolimits_{Q_{R}}\left\vert {Xv}\right\vert ^{2}dz}\leq\varepsilon
{\iint\nolimits_{Q_{4R}}\left\vert {Xu}\right\vert ^{2}}dz.
\]
Applying Lemma \ref{l3.6} to $w$,%

\begin{align*}
{\iint\nolimits_{Q_{\rho}}\left\vert {Xu}\right\vert ^{2}dz} &  \leq
2{\iint\nolimits_{Q_{\rho}}\left\vert {Xv}\right\vert ^{2}dz}+2{\iint
\nolimits_{Q_{\rho}}\left\vert {Xw}\right\vert ^{2}dz}\\
&  \leq c{\iint\nolimits_{Q_{\rho}}\left\vert {Xv}\right\vert ^{2}dz}%
+c(\frac{\rho}{R})^{Q+2}{\iint\nolimits_{Q_{R}}\left\vert {Xw}\right\vert
^{2}dz}\\
&  \leq c{\iint\nolimits_{Q_{\rho}}\left\vert {Xv}\right\vert ^{2}dz}%
+c(\frac{\rho}{R})^{Q+2}{\iint\nolimits_{Q_{R}}\left\vert {Xu}\right\vert
^{2}dz}\\
&  \leq c((\frac{\rho}{R})^{Q+2}+\varepsilon){\iint\nolimits_{Q_{R}}\left\vert
{Xu}\right\vert ^{2}dz}.
\end{align*}
The proof is reached by using Lemma 2.6.
\end{proof}

Next we discuss estimates of weak solutions of (\ref{1.1}) in parabolic cylinders.

\begin{theorem}
\label{T4.2}Under the assumption (H), let $u\in V_{2}(Q_{T})$ be a weak
solution of (\ref{1.1}) in $Q_{T}$ and $u=0$ on $\partial_{p}Q_{R}$. Suppose
that there exist $\lambda$ and $\gamma$ such that $\lambda<\gamma<Q+2$ and the
function $\frac{{r^{\gamma-\lambda}}}{{\varphi^{2}(r)}}$ is almost increasing.
Then $Xu\in L_{\varphi}^{2,\lambda}(Q_{\rho})$. Furthermore, for any $\rho\leq
R$, $Q_{R}\subset Q_{T}$, it follows
\[
{\iint\nolimits_{Q_{\rho}}\left\vert {Xu}\right\vert ^{2}dz}\leq c\frac
{{\rho^{\lambda}\varphi^{2}(\rho)}}{{R^{\lambda}\varphi^{2}(R)}}%
{\iint\nolimits_{Q_{R}}\left\vert {Xu}\right\vert ^{2}dz}+c\rho^{\lambda
}\varphi^{2}(\rho)(\left\Vert f\right\Vert _{L_{\varphi}^{2,\lambda}}%
^{2}+\left\Vert g\right\Vert _{L_{\varphi}^{2,\lambda}}^{2}).
\]

\end{theorem}

\begin{proof}
Let $w$ be a weak solution to the system%
\begin{equation}
\left\{
\begin{array}
[c]{c}%
w_{t}^{i}+X_{\alpha}^{\ast}(a_{ij}^{\alpha\beta}X_{\beta}w^{j})=0,\mathrm{{in}%
}\text{ }Q_{R},\\
w=u,\text{ \ \ \ \ \ \ \ }\mathrm{{on}}\text{ }\partial_{p}Q_{R},
\end{array}
\right.  \label{4.4}%
\end{equation}
Then $v=u-w$ satisfies
\begin{equation}
\left\{
\begin{array}
[c]{l}%
v_{t}^{i}+X_{\alpha}^{\ast}(a_{ij}^{\alpha\beta}X_{\beta}v^{j})=g_{i}%
+X_{\alpha}^{\ast}f_{i}^{\alpha},\mathrm{{in}}\text{ }Q_{R},\\
v=0,\text{ \ \ }\mathrm{{on}}\text{ }\partial_{p}Q_{R}.
\end{array}
\right.  \label{4.5}%
\end{equation}
Multiplying both sides of the system in (\ref{4.5}) by $v^{i}$ and integrating
on $Q_{R}$,
\[
{\iint\nolimits_{Q_{R}}(v_{t}^{i}v^{i}+a_{ij}^{\alpha\beta}X_{\beta}%
v^{j}X_{\alpha}v^{i})}dz={\iint\nolimits_{Q_{R}}(g_{i}v^{i}+f_{i}^{\alpha
}X_{\alpha}v^{i})dz}.
\]
Using (\ref{3.18}),%
\begin{align*}
&  {\iint\nolimits_{Q_{R}}(v_{t}^{i}v^{i}+a_{ij}^{\alpha\beta}X_{\beta}%
v^{j}X_{\alpha}v^{i})dz}\\
&  \leq C_{\varepsilon}{\iint\nolimits_{Q_{R}}(\left\vert g\right\vert
^{2}+\left\vert f\right\vert ^{2})dz}+\varepsilon{\iint\nolimits_{Q_{R}%
}\left\vert v\right\vert ^{2}dz}+\varepsilon{\iint\nolimits_{Q_{R}}\left\vert
{Xv}\right\vert ^{2}dz}\\
&  \leq2\varepsilon{\iint\nolimits_{Q_{R}}\left\vert {Xv}\right\vert ^{2}%
dz}+C_{\varepsilon}{\iint\nolimits_{Q_{R}}(\left\vert g\right\vert
^{2}+\left\vert f\right\vert ^{2})dz}.
\end{align*}
Since ${\iint\nolimits_{Q_{R}}v_{t}^{i}v^{i}dz}=\int_{B_{R}}{dx\int
_{t_{0}-R^{2}}^{t_{0}}{v^{i}dv^{i}}}\geq0$ and (\ref{1.2}), it yields
\begin{equation}
{\iint\nolimits_{Q_{R}}\left\vert {Xv}\right\vert ^{2}dz}\leq C_{\varepsilon
}{\iint\nolimits_{Q_{R}}(\left\vert g\right\vert ^{2}+\left\vert f\right\vert
^{2})dz}.\label{4.6}%
\end{equation}
Thanks to Theorem \ref{4.1}, we have%
\begin{align*}
{\iint\nolimits_{Q_{\rho}}\left\vert {Xu}\right\vert ^{2}dz} &  \leq
2{\iint\nolimits_{Q_{\rho}}\left\vert {Xw}\right\vert ^{2}dz}+2{\iint
\nolimits_{Q_{\rho}}\left\vert {Xv}\right\vert ^{2}dz}\\
&  \leq c(\frac{\rho}{R})^{\mu}{\iint\nolimits_{Q_{R}}\left\vert
{Xw}\right\vert ^{2}dz}+c{\iint\nolimits_{Q_{R}}\left\vert {Xv}\right\vert
^{2}dz}\\
&  \leq c(\frac{\rho}{R})^{\mu}{\iint\nolimits_{Q_{R}}\left\vert
{Xu}\right\vert ^{2}dz}+c{\iint\nolimits_{Q_{R}}\left\vert {Xv}\right\vert
^{2}dz}\\
&  \leq c(\frac{\rho}{R})^{\mu}{\iint\nolimits_{Q_{R}}\left\vert
{Xu}\right\vert ^{2}dz}+C_{\varepsilon}{\iint\nolimits_{Q_{R}}(\left\vert
g\right\vert ^{2}+\left\vert f\right\vert ^{2})dz}\\
&  \leq c(\frac{\rho}{R})^{\mu}{\iint\nolimits_{Q_{R}}\left\vert
{Xu}\right\vert ^{2}dz}+c\varphi^{2}(R)R^{\lambda}(\left\Vert f\right\Vert
_{L_{\varphi}^{2,\lambda}}^{2}+\left\Vert g\right\Vert _{L_{\varphi
}^{2,\lambda}}^{2}).
\end{align*}
Now letting $H(\rho)={\iint\nolimits_{Q_{\rho}}\left\vert {Xu}\right\vert
^{2}dz}$, $H(R)={\iint\nolimits_{Q_{R}}\left\vert {Xu}\right\vert ^{2}dz}$,
$B=\left\Vert f\right\Vert _{L_{\varphi}^{2,\lambda}}^{2}+\left\Vert
g\right\Vert _{L_{\varphi}^{2,\lambda}}^{2}$, $F(R)=\varphi^{2}(R)R^{\lambda}$
and $\beta=Q+2,$ and noting that the function $\frac{{\rho^{\gamma}}}%
{{F(\rho)}}$ is almost increasing in $(0,R_{0}]$, we have by Lemma \ref{l2.1}
that
\[
{\iint\nolimits_{Q_{\rho}}\left\vert {Xu}\right\vert ^{2}dz}\leq c\frac
{{\rho^{\lambda}\varphi^{2}(\rho)}}{{R^{\lambda}\varphi^{2}(R)}}%
{\iint\nolimits_{Q_{R}}\left\vert {Xu}\right\vert ^{2}dz}+c\rho^{\lambda
}\varphi^{2}(\rho)(\left\Vert f\right\Vert _{L_{\varphi}^{2,\lambda}}%
^{2}+\left\Vert g\right\Vert _{L_{\varphi}^{2,\lambda}}^{2}).
\]
This proof is completed.
\end{proof}

\begin{proof}
[\textit{Proof for Theorem 1.1}]By Theorem \ref{T4.2} and the cutoff function
technique, it is easy to see that Theorem \ref{T1.1} is true, and we omit the details.
\end{proof}

\mbox{} Department of Applied mathematics, Northwestern polytechnical
university; Key laboratory of Space Applied physics and chemistry, Ministry of
Education, Xi'an, shaanxi, 710129, China\newline E-mail addresses:
dongyan8617@sina.com\newline E-mail addresses: zhumaochun2006@126.com\newline
E-mail addresses: pengchengniu@nwpu.edu.cn

\end{document}